\documentclass[10pt]{amsart}
\usepackage{amsmath}
\usepackage{enumerate}
\usepackage[colorlinks=true]{hyperref}
\newtheorem{theorem}{Theorem}[section]

\newtheorem{lemma}[theorem]{Lemma}
\newtheorem{proposition}[theorem]{Proposition}
\newtheorem{example}[theorem]{Example}
\newtheorem{definition}[theorem]{Definition}
\newtheorem{remark}{Remark}
\numberwithin{equation}{section}

\newcommand{\abs}[1]{\left\vert#1\right\vert}

\newcommand{\real}{\mathbb R}

\def\natu{\mathbb N}

\def\D{\nabla}

\newcommand{\sgn}[1]{{\rm sgn}(#1)}
\def\D{\nabla}

\def\D{\nabla}
\newcommand{\RR}{\mathbb{R}}

\newcommand{\Q}{\mathcal{Q}}

\begin{document}
	
	\title[Some results for Semi-stable solutions of $k$-Hessian equations]
	{Some results for Semi-stable radial solutions of $k$-Hessian equations}
	\author{Miguel Angel Navarro and Justino S\'anchez}
	\address{Departamento de Estat\'{i}stica, An\'alise Matem\'atica e Optimizaci\'on\\Universidade de Santiago de Compostela\\Santiago de Compostela 15782, Spain.}
	\address{Departamento de Matem\'{a}ticas\\
		Universidad de La Serena\\
		Avenida Cisternas 1200, La Serena, Chile.}	
	\email{miguel.navarro.burgos@gmail.com \and jsanchez@userena.cl}

	\begin{abstract}
		We devote this paper to study semi-stable nonconstant radial solutions of $S_k(D^2u)=w(\abs{x})g(u)$ on the Euclidean space $\real^n$. We establish pointwise estimates and necessary conditions for the existence of such solutions (not necessarily bounded) for this equation. For bounded solutions we estimate their asymptotic behavior at infinity. All the estimates are given in terms of the spatial dimension $n$, the values of $k$ and the behavior at infinity of the growth rate function of $w$. 
	\end{abstract}
	
	\maketitle
	\section{Introduction and main results}
	This paper deals with the semi-stability of nonconstant radial
	solutions of
	\begin{equation}\label{mainequation}
		S_k(D^2u)=w(\abs{x})g(u),
	\end{equation}
	posed in $\real^n$, where $n\geq 1$, $k\in\{1,\ldots,n\}$, the function $g\in C^1(\real)$ is nonnegative and nonincreasing, and $w$ is a nonnegative radial function that satisfies some suitable conditions. The operator $S_k(D^2u)$ is the $k$-Hessian of $u$, which is defined by the sum of all $k$-th principal minors of the Hessian matrix. Alternatively, $S_k(D^2u)$ is the $k$-th elementary symmetric polinomial of the eigenvalues of the Hessian matrix $D^2 u$. 
	According to \cite{CaNS85}, to ensure ellipticity of equation \eqref{mainequation}, we consider $k$-admissible or $k$-convex solutions, i.e., functions that belongs to 
	\begin{equation*}
		\Phi^k:=\left\lbrace u\in C^2(\real^n): S_l(D^2 u)\geq 0,\, l=1,2,\ldots,k\right\rbrace.
	\end{equation*}
	
	We point out that there are no previous works concerning semi-stable solutions to the $k$-Hessian operator in the whole space $\real^n$ for general nonlinearities and $k\neq 1$. In fact, to the best of our knowledge, the only results devoted to stable solutions to \eqref{mainequation} are contained in the work \cite{WaLe19} and they concern only the particular cases $w(r)\equiv 1$ and $g(u)$ a power nonlinearity. More precisely, in \cite{WaLe19}, the authors gave a definition of stable radial solutions of the $k$-Hessian equation $F_{k}(D^2V)=(-V)^p$ in $\RR^n$, where $F_{k}(D^2V)=S_k(D^2V)$. They stablished connections between stability and certain critical exponents of Joseph-Lundgren type available for $k$-Hessian operators. Their tools also include Wolff potentials.
	
	For existence and non-existence results for equations of the form \eqref{mainequation} we refer to \cite{ClMM98}, where a special emphasis was put in the model equation 
	\begin{equation*}
	(-1)^kH_k u=a(\abs{x})\abs{u}^{q-1}u,
	\end{equation*}
	being
	\begin{equation*}
		H_k=\frac{k}{\binom{n-1}{k-1}}T_k,\mbox{ for }1\leq k \leq n\mbox{ (integer)}\mbox{ and }T_k=S_k(\nabla^2).
	\end{equation*}
	
	Some existence and non-existence results for radial solutions are given in terms of an integral condition involving the function $a$. In \cite{Lei15}, among others results, the author construct explicit negative solutions of the equation $F_{k}(D^2V)=R(x)(-V)^q$ in $\RR^n$, where $F_{k}(D^2V)=S_k(D^2V)$ and $R(x)$ is a radial function that satisfies $C^{-1}\leq R(x)\leq C$ for some constant $C>1$. See \cite[Theorem 4.2]{Lei15}.
	
	Throughout this work, we identify a radial solution $u$ by their one variable representant, that is, $u(x) = u(r),\, \abs{x}=r$. At the point $x=(r,0,...,0)$, the eigenvalues of $D^2 u$ are $\lambda_1=u''$, which is simple, and $\lambda_2=\frac{u'}{r}$, which has multiplicity $n-1$, where by abuse of notation, we write $u'$ or $\partial_r u$ as the radial derivative of a radial function $u$. Thus the $k$-Hessian operator acting on radially symmetric $C^2$ functions can be written as 
	\begin{equation}\label{eq:radk-Hess}
		S_k(D^2 u)=c_{n,k}\lambda_2^{k-1}\left(n\lambda_2+k\left(\lambda_1-\lambda_2\right)\right)=c_{n,k}r^{1-n}\partial_r\left (r^n\lambda_2^k\right ),\, r>0,
	\end{equation}
	where $c_{n,k}$ is defined by $c_{n,k}=\binom{n}{k}/n$.  
	
	\begin{remark}\label{nondecre}
		Note that, if $u$ is a radial solution of \eqref{mainequation} then, in particular, $S_1(D^2 u)=r^{1-n}(r^{n-1}u')'\geq 0$. Thus $G(r)=r^{n-1}u'$
		is nondecreasing, since $G(0)=0$, we deduce that $G\geq 0$ and hence $u$ is nondecreasing. As a consequence, $\lambda_2\geq 0$.
	\end{remark}

	\begin{definition} We say that a radial solution $u\in\Phi^k$ of \eqref{mainequation} is semi-stable if		
	\begin{equation}\label{semistable}
		\Q_u(\xi):=\int_{\real^n}{kc_{n,k}\abs{x}^{1-k}\abs{\nabla u}^{k-1}\abs{\nabla\xi}^2+w(\abs{x})g'(u)\xi^2}\geq 0,
	\end{equation}
	for every radially symmetric function $\xi\in C_c^{1}(\real^n)$.
	\end{definition}

	In this paper, we establish pointwise estimates and necessary conditions for the existence of semi-stable solutions (not necessarily bounded) of \eqref{mainequation}. For bounded solutions we estimate their asymptotic behavior at infinity. All the estimates are given in terms of the spatial dimension $n$, the values of $k$ and the behavior at infinity of the growth rate function of $w$, that is, the function $W(r):=\frac{rw'(r)}{w(r)}$. 
	
	We now stablish our precise assumptions on the weight $w(r)$ that we will assume throughout the paper:	
	\begin{equation}\label{ineq:H1rwrw}
		\begin{cases}
			\mbox{The functions } w \mbox{ and } W \mbox{ belongs to }C^1\left(\real^n\setminus\{0\}\right)
			\mbox{ and } W(r) \mbox{ is}\\
			\mbox{nonincreasing on }(0,\infty). 
			\mbox{ Set }\Gamma=\lim\limits_{r\to 0}W(r),\,\gamma=\lim\limits_{r\to +\infty}W(r)\\ 
			\mbox{and assume that }\Gamma \mbox{ and }\gamma \mbox{ are finite}.
		\end{cases}
	\end{equation}

In order to state our main results, we need the following notation:

Let
	\begin{equation}\label{eq:defalpha}
			\alpha(r)=\frac{\int_1^r{\left\lbrace \frac{v(s)+\sqrt{v(s)^2+\frac{n}{k}v(s)+\frac{n}{k}-1 -sv'(s)}}{s}\right\rbrace ds}}{\log r},
		\end{equation}
%
		\begin{equation}\label{eq:deltank}
			\delta(r)=\frac{-n+2\alpha(r)+2k+2}{k+1},\,\forall r\geq 1\mbox{ and }		\delta_{\infty}(\gamma)=\lim\limits_{r\to +\infty}{\delta(r)},
		\end{equation}	
		where
\begin{equation}\label{eq:defv}
v(s)=\frac{k-1+W(s)}{k+1},\,\forall s\geq 0.
\end{equation}	
		
Note that when $k=1$ and $w\equiv 1$, $\delta(r)$ take the constant value $-\frac{n}{2}+\sqrt{n-1}+2$, which plays a crucial role in pointwise estimates. We refer to \cite{MR2355457} by Villegas for details. A great difference here is that $\delta(r)$ is a variable exponent, this produce additional technical difficulties.  
		
	Our main results are
	\begin{proposition}\label{gradnonulo}
		Let $n>k(k+1)/(2k+\gamma)$ and $u$ be a semi-stable nonconstant radial solution of \eqref{mainequation}. 
		Then $\abs{\D u(x)}>0$ for all $\abs{x}>0$.
	\end{proposition}
	
	We have the following pointwise estimate for not necessarily bounded solutions of \eqref{mainequation}.
	\begin{theorem}\label{nbgrowth}
		Let $w\in C^1\left(\real^n\setminus\{0\}\right)$ be a radial function that satisfies \eqref{ineq:H1rwrw}, $n>k(k+1)/(2k+\gamma)$ and $u$ be a semi-stable nonconstant radial solution of \eqref{mainequation} (not necessarily bounded). Then, there
		exist $M>0$ and $r_0\geq 1$ depending on $u$ and $w$ such that for all $r\geq r_0$,
		\begin{equation}\label{ingrowth}
			\abs{u(r)}\geq M\begin{cases}
				r^{\delta(r)}&\mbox{if }n\neq 2\left(k+\frac{2\gamma}{k}+4\right),\\
				\log r&\mbox{if }n=2\left(k+\frac{2\gamma}{k}+4\right).
			\end{cases}
		\end{equation}
	\end{theorem}
	
	\begin{remark}\label{rmk1}
	
		This theorem is sharp for some $\beta\in\RR,\, \mu=\frac{\sigma}{k}+2$ and 
		\begin{equation*}
			u_{\beta} (r):=\begin{cases}
				\sgn{\beta}\left(1+r^{\mu}\right )^{\frac{\beta}{\mu}} &\mbox{if } \beta
				\neq 0,\forall r\geq 0,\\\\
				\left(\frac{1}{\mu}\right)\log{\left (1+r^{\mu}\right )}&\mbox{if } \beta=0,\forall r\geq 0,
			\end{cases}
		\end{equation*}
	 where $\sigma>-k$ and $w(r)=r^{\sigma}$, as shown in the Appendix.
		
	\end{remark}
	
	\begin{theorem}\label{bogrowth}
		Let $w\in C^1\left(\real^n\setminus\{0\}\right)$ be a radial function that satisfies \eqref{ineq:H1rwrw}, $n\geq 2$ and $u$ be a semi-stable nonconstant bounded radial solution of \eqref{mainequation}. Then,
		
		\begin{enumerate}[$i)$]
			\item $n>2\left(k+\frac{2\gamma}{k}+4\right)$.
			\item There exists $u_\infty=\lim\limits_{r\rightarrow +\infty}{u(r)}\in\RR$ and $M>0$ depending on $u$ and $w$ such that for all $r\geq 1$,
			\begin{equation}\label{growth}
				\abs{u(r)-u_\infty} \geq M r^{\delta(r)}.
			\end{equation}		
		\end{enumerate}
	\end{theorem}
	\begin{remark}\label{rmk2}
		Theorem \ref{bogrowth} is sharp (see Example \ref{exam:family} in the Appendix).	
	\end{remark}
	
	\begin{remark}
		See \cite{MR2355457} to compare our results with the semilinear case ($k=1$) and for related equations involving the $p$-{L}aplacian operator, see \cite{MR2552149, Fazly12, paper3}.
	\end{remark}
	
	This paper is organized as follows. In Section 2 we prove our main results, Proposition \ref{gradnonulo}, Theorem \ref{nbgrowth} and Theorem \ref{bogrowth}. We conclude the paper by presenting in the Appendix some examples of functions $w$ and $g$ for which our theorems are sharp.
	
	\section{Proof of the main results}
	We claim that if $u$ is a $C^2$ radial solution, then $u'\in C_{loc}^{0,1}(\RR^n)$. To prove our claim, we first observe that $u'\in C^0(\RR^n)\cap C^1(\RR^n\setminus\{0\})$ with $u'(0)=0$ and $\abs{u'(x)}\leq C_R\abs{x}$ in any open ball $B_R$ of radius $R>0$. Now
	from \eqref{mainequation} and \eqref{eq:radk-Hess} we obtain, for $r>0$,
	\begin{equation*}
		kc_{n,k}\left(\frac{u'}{r}\right)^{k-1}\left(\left(\frac{n-k}{k}\right)\left(\frac{u'}{r}\right)+u''\right)=wg(u),
	\end{equation*}
	from which we deduce that
	\begin{equation*}\label{eq:d2u}
		u''=\left(\frac{wg(u)}{kc_{n,k}}\right) \left(\frac{u'}{r}\right)^{1-k}-\left(\frac{n-k}{k}\right)\left(\frac{u'}{r}\right),\,r>0.
	\end{equation*}
	
	Thus, for any $R>0$, the function $u''$ is bounded on $B_R\setminus\{0\}$. Also for any $x,y\in\RR^n$ such that $R>\abs{y}>\abs{x}>0$, we have
	\begin{equation*}
		\begin{split}
			\abs{u'(\abs{y})-u'(\abs{x})}&\leq\int_{\abs{x}}^{\abs{y}}\abs{u''(t)}dt\\
			&\leq\sup_{\xi\in [\abs{x},\abs{y}]}\abs{u''(\xi)}(\abs{y}-\abs{x})\\
			&\leq\sup_{z\in B_{R}\setminus\{0\}}\abs{u''(\abs{z})}\abs{y-x}.
		\end{split}
	\end{equation*}
	
	We conclude that $u'\in C_{loc}^{0,1}(\RR^n)$. This will be used in the proof of Lemma \ref{essential0} below.	
	\begin{lemma}\label{essential0}
		Let $n\geq 2$ and $u$ be any radial solution of
		\eqref{mainequation}. Then 		
		\begin{equation}\label{eq:des_semi_stable0}
			\mathcal{Q}_u(u'\eta)=kc_{n,k}\int_{\real^n}{\lambda_2^{k+1}\left\lbrace\abs{\abs{x}\D\eta+v\frac{x\eta}{\abs{x}}}^2-\left(v^2+\frac{n}{k}v+\frac{n}{k}-1-\left(x,\D v\right)\right)\eta^2\right\rbrace },
		\end{equation}
		for every radially symmetric function $\eta\in (H_c^1\cap L_{loc}^{\infty})\left(\RR^n\right)$, where $v$ is defined as in \eqref{eq:defv}. Here $(\cdot, \cdot)$ denotes the standard scalar product in $\RR^n$.
	\end{lemma}
	
	\begin{proof}
		Let $\eta\in H_c^1(\RR^n)\cap L_{loc}^\infty(\RR^n)$ and $\zeta\in C_{loc}^{0,1}(\RR^n)$ be radial functions. Then, by a standard density argument, we can take $\xi=\zeta\eta\in H_c^1(\RR^n)\cap L_{loc}^\infty(\RR^n)$ in \eqref{semistable} to obtain
		\begin{equation}\label{eq:Qzeta}
			\Q_u(\zeta\eta)=\int_{\RR^n}{kc_{n,k}\abs{\frac{u'}{\abs{x}}}^{k-1}\abs{\D(\zeta\eta)}^2+w(\abs{x})g'(u)(\zeta\eta)^2}.
		\end{equation}
		
		Thus, as $u$ and $\zeta$ are radial functions, we get 
		\begin{equation}\label{eq:radialzetaeta}
			\abs{\frac{u'}{\abs{x}}}^{k-1}\abs{\D(\zeta\eta)}^2=\lambda_2^{k-1}\left(\zeta^2\abs{\nabla\eta}^2+\left(2\zeta\eta(\nabla\zeta,\nabla\eta)+\eta^2\abs{\nabla\zeta}^2\right)\right),
		\end{equation}
		where $\lambda_2=u'/{\abs{x}}.$
		
		From \eqref{eq:radk-Hess} and differentiating \eqref{mainequation} with respect to $r$, we obtain
		\begin{equation*}
			\begin{split}
				c_{n,k}^{-1}wg'(u)u'&=-c_{n,k}^{-1}w'g(u)-\frac{n-1}{r^2}\left(\frac{u'^k}{r^{k-1}}\right)+\Delta\left(\frac{u'^k}{r^{k-1}}\right)\\
				&=-c_{n,k}^{-1}wg(u)\frac{W}{r}-\frac{(n-1)\lambda_2^k}{r}+\Delta\left(\lambda_2^{k-1}u'\right),\,r>0.
			\end{split}
		\end{equation*}
		
		Then, multiplying the latter equation by $u'\eta^2$, integrating by parts and using $\zeta=u'$ in \eqref{eq:Qzeta}, \eqref{eq:radialzetaeta}, and taking into account that $\lambda_2\in L_{loc}^\infty
		(\RR^n)$, we have		
		\begin{equation*}\label{lemmaig1}
			\begin{split}
				\frac{1}{kc_{n,k}}\int_{\real^n}{wg'(u)(u'\eta)^2}&=-\int_{\real^n}{\lambda_2^{k-1}\abs{\D\left(u'\eta\right)}^2}\\
				&+\int_{\real^n}{\lambda_2^{k+1}\left\lbrace\abs{x}^2\abs{\D\eta}^2+\left(\frac{k-1}{k+1}\right)\left(x,\D\eta^2\right)\right\rbrace }\\
				&-\int_{\real^n}{\lambda_2^{k+1}\left\lbrace \frac{2n-k-1}{k+1}+\frac{n}{k}W\right\rbrace \eta^2}\\
				&-\frac{1}{k+1}\int_{\real^n}{W\left(x,\D\lambda_2^{k+1}\right)\eta^2}\\
				&=-\int_{\real^n}{\lambda_2^{k-1}\abs{\D\left(u'\eta\right)}^2}\\
				&+\int_{\real^n}{\lambda_2^{k+1}\left\lbrace\abs{x}^2\abs{\D\eta}^2+\left(\frac{k-1+W}{k+1}\right)\left(x,\D\eta^2\right)\right\rbrace }\\
				&-\int_{\real^n}{\lambda_2^{k+1}\left\lbrace \frac{2n-k-1+\frac{n}{k}W-\left(x,\D W\right)}{k+1}\right\rbrace \eta^2}.
			\end{split}
		\end{equation*}
			
		Now from \eqref{eq:defv}, we have
		\begin{equation*}\label{eq:essrdww}
			W=(k+1)v+1-k\mbox{ and }\left(x,\D W\right)=(k+1)\left(x,\D v \right).
		\end{equation*}
		
		Hence
		\begin{equation*}\label{lemmaig1}
			\begin{split}
				\int_{\real^n}{\lambda_2^{k-1}\abs{\D\left(u'\eta\right)}^2}+\frac{1}{kc_{n,k}}\int_{\real^n}{wg'(u)(u'\eta)^2}&=\int_{\real^n}{\lambda_2^{k+1}\left\lbrace\abs{x}^2\abs{\D\eta}^2+v\left(x,\D\eta^2\right)\right\rbrace }\\
				&-\int_{\real^n}{\lambda_2^{k+1}\left\lbrace\frac{n}{k}v+\frac{n}{k}-1-\left(x,\D v\right)\right\rbrace \eta^2},
			\end{split}
		\end{equation*}
		which concludes the proof.
		
	\end{proof}	
	
	Multiplying the equation \eqref{eq:radk-Hess} by $r^{n+\frac{n}{k}-1}\lambda_2$, we obtain
	\begin{equation*}
		\left(\left(r^{\frac{n}{k}}\lambda_2\right)^{k+1}\right)'=(k+1)(kc_{n,k})^{-1}wg(u)r^{n+\frac{n}{k}-1}\lambda_2.
	\end{equation*}
	
	Since $r\lambda_2=u'$, the above equation is equivalent to
	\begin{equation}\label{conveform}
		\left(\left(r^{\frac{n}{k}-1}u'\right)^{k+1}\right)'=(k+1)(kc_{n,k})^{-1}r^{n-1}wg(u)r^{\frac{n}{k}-1}u'.
	\end{equation}
	
	On the other hand, from \eqref{ineq:H1rwrw} and \eqref{eq:defv}, we have
	\begin{equation}\label{eq:defvinf}
		v(r)\geq v_{\infty}:=\lim_{r\to +\infty}{v(r)}=\frac{k-1+\gamma}{k+1},\,\forall r\geq 0.
	\end{equation}
	
	Thus, if $n>k(k+1)/(2k+\gamma)$, it follows that	
	\begin{equation}\label{ineq:nkvtvtc1}
		\frac{n}{k}v(r)+\frac{n}{k}-1\geq\frac{n}{k}\left(v_{\infty}+1\right)-1=\frac{n(2k+\gamma)}{k(k+1)}-1>0,\,\forall r\geq 0.
	\end{equation}
	
	\begin{proof}[Proof of Proposition \ref{gradnonulo}] We follows an argument similar to that of Proposition 1 in \cite{FaNa20}. Let $n>k(k+1)/(2k+\gamma)$ and let $u$ be a semi-stable nonconstant radial solution of
		\eqref{mainequation}. Arguing by contradiction, assume that $u'(r_0) = 0$ for some $r_0 > 0$, so $u'\chi_{B(0,r_0)}\in H_c^1(\RR^n)\cap L_{loc}^\infty(\RR^n)
		$ and from \eqref{eq:des_semi_stable0}, we get that
		\begin{equation*}
			\mathcal{Q}_u\left(u'\chi_{B(0,r_0)}\right)=-kc_{n,k}\int_{B(0,r_0)}\left(\frac{n}{k}v+\frac{n}{k}-1-\left(x,\D v\right)\right)r^2\lambda_2^{k+3}\leq 0,
		\end{equation*}
		where we have used the fact that $\lambda_2\geq 0$, \eqref{ineq:nkvtvtc1} and the monotonicity of $v$ to get the above inequality. 
		Thus, the semi-stability of $u$ implies that $u(r)=u_0$ for all $r\in [0,r_0]$.
		
		Let $v(r)=(r^{\frac{n}{k}-1}u')^{k+1}$. From \eqref{conveform}, we have the following problem
		\begin{equation*}
			\left\lbrace\begin{array}{rcll}
			u'&=&r^{1-\frac{n}{k}}v^{\frac{1}{k+1}}& \quad r>r_0>0,\\
			v'&=&(k+1)(kc_{n,k})^{-1}r^{n-1}wg(u)v^{\frac{1}{k+1}}& \quad r>r_0>0,\\
			u(r_0)&=&u_0,&\\
			v(r_0)&=&0.&\\
			\end{array}\right.			
		\end{equation*}
		
		Finally, by Cauchy's theorem, we get that $u=u_0$ for any $r>0$, a contradiction.
	\end{proof}

	We adapt some estimates given by Villegas in \cite{MR2355457} for the nonweighted semilinear equation. Here is one of our main integral estimates.
	\begin{lemma}\label{essentialg}
		Let $w\in C^1\left(\real^n\setminus\{0\}\right)$ be a radial function that satisfies \eqref{ineq:H1rwrw}, $n>k(k+1)/(2k+\gamma)$ and $u$ be a semi-stable nonconstant radial solution of
		\eqref{mainequation}. Then, there exists $K>0$ depending on $u$ and $w$ and such that		
		\begin{equation}\label{ineq:inequality}
			\frac{}{}
			\int_r^R{\frac{ds}{s^{n-k}({u'(s)})^{k+1}}}\leq Kr^{-2\alpha(r)}\,\ \ \forall R>r\geq 1.
		\end{equation}		
	\end{lemma}
	
	\begin{proof}
	
		From Proposition \ref{gradnonulo}, we have that $u'(r)\neq 0$ for all $r>0$. Furthermore, $u'>0$ on $(0,\infty)$ (see Remark \ref{nondecre}).
		
		Next we show that $\alpha(r)$ defined in \eqref{eq:defalpha} is strictly positive for all $r\geq 1$. For this, let $\kappa\geq 0$ and define the function
		\begin{equation*}
			\Psi(r,\kappa):=v(r)+\sqrt{v(r)^2+\frac{n}{k}v(r)+\frac{n}{k}-1 +\kappa}.
		\end{equation*}
		
		From \eqref{ineq:nkvtvtc1}, we have
		\begin{equation}\label{ineq:vn2k}
			v(r)+\frac{n}{2k}>\frac{k}{n}\left(1-\frac{n}{k}\right)+\frac{n}{2k}=\frac{k^2+(n-k)^2}{2kn}>0. 
		\end{equation}
		
		Then
		\begin{equation}\label{eq:dPhi}
			\frac{\partial\Psi}{\partial r}(r,\kappa)=\left(1+\frac{v(r)+\frac{n}{2k}}{\sqrt{v(r)^2+\frac{n}{k}v(r)+\frac{n}{k}-1 +\kappa}}\right)v'(r)\leq 0\mbox{ for }r>0.
		\end{equation}
				
		Thus, from \eqref{ineq:H1rwrw}, \eqref{ineq:nkvtvtc1}, \eqref{ineq:vn2k} and \eqref{eq:dPhi}, we obtain
		\begin{equation*}
			\Psi(r,-rv'(r))>v(r)+\abs{v(r)}\geq 0,\,\forall r\geq 0,
		\end{equation*}
		and 
		\begin{equation*}
			\begin{split}
				\Psi(r,-rv'(r))&\geq\Psi(r,0)=v(r)+\sqrt{\left(\frac{v(r)+1}{k}\right)\left(k\left(v(r)-1\right)+n\right)}\\
				&\geq\lim_{r\to +\infty}{\Psi(r,0)}=v_{\infty}+\sqrt{\frac{n}{k}(v_{\infty}+1)+(v_{\infty}+1)(v_{\infty}-1)}\\
				&>\frac{k-1+\gamma}{k+1}+\sqrt{1+\frac{(2k+\gamma)(\gamma-2)}{(k+1)^2}}\\
				&=\frac{k-1+\gamma+\abs{k-1+\gamma}}{k+1}\geq 0.
			\end{split}
		\end{equation*}
		
		Therefore by the previous inequalities, we get
		\begin{equation}\label{ineq:Malpm}
			\begin{split}
				\alpha(r)&=\frac{\int_1^r{\left\lbrace \frac{v(s)+\sqrt{v(s)^2+\frac{n}{k}v(s)+\frac{n}{k}-1 -sv'(s)}}{s}\right\rbrace ds}}{\log r}\\
				&\geq v_{\infty}+\sqrt{\frac{n}{k}(v_{\infty}+1)+(v_{\infty}+1)(v_{\infty}-1)}>0\mbox{ for }r\geq 1.
			\end{split}
		\end{equation}
		
		We now fix $R>r\geq 1$ and consider the function
		\begin{equation*}
			\eta (t)=\begin{cases}
				1 & \mbox{if } 0\leq t \leq 1, \\
				t^{-\alpha(t)} & \mbox{if } 1<t\leq r,\\
				\frac{r^{-\alpha(r)}\int_t^R{\frac{ds}{s^{n-k}(u'(s))^{k+1}}}}{\int_r^R{\frac{ds}{s^{n-k}(u'(s))^{k+1}}}}&
				\mbox{if } r<t\leq R,\\
				0 &\mbox{if } R<t<\infty.
			\end{cases}
		\end{equation*}
		
		Since, $u$ is a semi-stable solution, from \eqref{eq:des_semi_stable0}, we have 		
		\begin{equation}\label{ineq:I123}
			\begin{split}
				0\leq (\omega_n kc_{n,k})^{-1}\mathcal{Q}_u(u'\eta)&= I_1+I_2+I_3\\
				&=\left(\int_0^1+\int_1^r+\int_r^R\right){\left\lbrace  t^{n-1}\lambda_2^{k+1}\left( \left(t\eta'\right)^2+2v(t)\left(t\eta\eta'\right)\right.\right. }\\ &-{\left.\left. \left(\frac{n}{k}v(t)+\frac{n}{k}-1-tv'(t)\right)\eta^2\right)dt\right\rbrace },
			\end{split}
		\end{equation}
		where $\omega_n$ is the measure area of the $n-1$ dimensional unit sphere $S^{n-1}$.
		
		Then by \eqref{ineq:nkvtvtc1}, we have
		\begin{equation}\label{ineq:I1}
			\begin{split}
				I_1&=-\int_0^1{t^{n-1}\lambda_2^{k+1}\left(\left(\frac{n}{k}v(t)+\frac{n}{k}-1-tv'(t)\right)\right)dt}\\
				&\leq-\left(\frac{n}{k}\left(v_{\infty}+1\right)-1\right)\int_0^1{t^{n-1}\lambda_2^{k+1}\,dt},
			\end{split}
		\end{equation}		
		and by \eqref{eq:defalpha}, we obtain
		\begin{equation}\label{eq:I2}
			\begin{split}
				I_2&=\int_1^r{t^{n-2\alpha(t)-1}\lambda_2^{k+1}\left( \left(t\left(\alpha(t)\log t\right)'-v(t)\right)^2\right.}\\
				&-{\left. \left(v(t)^2+\frac{n}{k}v(t)+\frac{n}{k}-1-tv'(t)\right)\right)dt }\\
				&=0.
			\end{split}
		\end{equation}
		
		Now to estimate $I_3$, we need to consider two cases according to the sign of $v$ in the interval $[r,R]$. For this, we rewrite $I_3$ as
		\begin{equation}\label{eq:I3+-}
			\begin{split}
				I_3&=I_{3}^{+}+I_{3}^{-}=\int_{\{t\in[r,R]:v(t)\geq0\}}+\int_{\{t\in[r,R]:v(t)<0\}},
			\end{split}
		\end{equation} 
		and we define $\Phi_R: [r,R]\rightarrow\mathbb{R}$ as	
		\begin{equation*}
			\Phi_R(t)=\int_t^R{s^{-(n+1)}\lambda_2^{-(k+1)}\,ds},
		\end{equation*}
		and recall that $\lambda_2=u'(s)/s$ for $s>0$.
		
		Thus for $t\in[r,R]$:
		\begin{itemize}
			\item If $v(t)\geq0$, then by \eqref{ineq:nkvtvtc1}, it follows that
			\begin{equation}\label{ineq:I3}
				I_3^{+}\leq \frac{r^{-2\alpha(r)}}{\Phi_R(r)}.
			\end{equation}
			\item If $v(t)<0$, then for all $\epsilon>0$ 
			\begin{equation*}
				-2v(t)\Phi_R(t)t^{-n}\lambda_2^{-(k+1)}\leq \epsilon^2v(t)^2\Phi_R(t)^2+\frac{t^{-2n}\lambda_2^{-2(k+1)}}{\epsilon^2},
			\end{equation*}
			and
			\begin{equation*}
				\begin{split}
					I_3^{-}&\leq \left(1+\frac1{\epsilon^2}\right) \frac{r^{-2\alpha(r)}}{\Phi_R(r)}-\frac{r^{-2\alpha(r)}}{\Phi_R(r)^2}\times\\
					&\times\int_{\{t\in[r,R]:v(t)<0\}}{t^{n-1}\lambda_2^{k+1}\left( -\epsilon^2v(t)^2+\frac{n}{k}v(t)+\frac{n}{k}-1-tv'(t)\right)\Phi_R(t)^2dt }.
				\end{split}
			\end{equation*}
			
			Moreover, by \eqref{ineq:H1rwrw} and \eqref{eq:defv}, we have that 
			\begin{equation*}
				\abs{v_{\infty}}>\abs{v(t)}>0\mbox{ for any }t\in\{t\in[r,R]:v(t)<0\}.
			\end{equation*}
			
			By \eqref{ineq:nkvtvtc1} we can pick $\epsilon\in\left(0,\sqrt{\frac{n}{k}v_{\infty}+\frac{n}{k}-1}/\abs{v_{\infty}}\right)$. Then applying the Young's inequality with $\varepsilon$ to $-v(t)\Phi_R(t) t^{-n}\lambda_2^{-(k+1)}$ with exponents 2 and 2, and using again \eqref{ineq:nkvtvtc1}, we have
			\begin{equation*}
				\begin{split}
					-\epsilon^2v(t)^2+\frac{n}{k}v(t)+\frac{n}{k}-1-tv'(t)&\geq-\epsilon^2v_{\infty}^2+\frac{n}{k}v_{\infty}+\frac{n}{k}-1>0, 
				\end{split}
			\end{equation*}
			and therefore
			\begin{equation}\label{ineq:CI3}
				I_3^{-}\leq\left(1+\frac{1}{\epsilon^2}\right) \frac{r^{-2\alpha(r)}}{\Phi_R(r)}.
			\end{equation}	  	
		\end{itemize}
		
		Finally, from \eqref{ineq:I123}-\eqref{eq:I2}, \eqref{ineq:I3} and \eqref{ineq:CI3}, the lemma follows.
	\end{proof}

	Applying Lemma \ref{essentialg} enables us to prove the following pointwise estimate.
	\begin{proposition}\label{prorand2r}
		
		Let $w\in C^1\left(\real^n\setminus\{0\}\right)$ be a radial function that satisfies \eqref{ineq:H1rwrw}, $n>k(k+1)/(2k+\gamma)$ and $u$ be a semi-stable nonconstant radial solution of
		\eqref{mainequation}. Then, there exists $K'>0$ depending on $u$ and $w$ and such that		
		\begin{equation}\label{rand2r}
			\abs{u(2r)-u(r)}\geq K'r^{\delta(r)},\quad\forall r\geq 1.
		\end{equation}		
	\end{proposition}
	
	\begin{proof}	
		Fix $r\geq 1$. Applying H\"older's inequality, Lemma \ref{essentialg} with $R=2r$ and recalling that $u'$ does not vanish in $(0,\infty)$, we deduce		
		\begin{equation*}
			\begin{gathered}
				\left (\int_1^{2}{ t^{-\frac{n-k}{k+2}}\,dt}\right )r^{-\frac{n-2k-2}{k+2}}=\int_r^{2r}{ t^{-\frac{n-k}{k+2}}\,dt}\\
				\leq \left (\int_r^{2r}{ \frac{dt}{t^{n-k}({u'(t)})^{k+1}}}\right )^{\frac1{k+2}} \left (\int_r^{2r}{{u'(t)}\,dt}\right )^{\frac{k+1}{k+2}}\\
				\leq \left(Kr^{-2\alpha(r)} \abs{u(2r)-u(r)}^{k+1}\right)^\frac{1}{k+2},
			\end{gathered}
		\end{equation*}
		which gives \eqref{rand2r}.
	\end{proof}
	
	\begin{proof}[Proof of Theorem \ref{nbgrowth}]
		
		Let $\delta(r)$ and $\delta_\infty(\gamma)$ as in \eqref{eq:deltank}. By \eqref{eq:defvinf} and \eqref{ineq:nkvtvtc1}, we obtain
		\begin{equation*}
			v_{\infty}+1=\frac{2k+\gamma}{k+1}>0\mbox{ and }k(v_{\infty}-1)+n>\frac{(n-k)^2}{n}\geq 0.
		\end{equation*}
		
		Thus, from \eqref{ineq:H1rwrw}, \eqref{eq:defalpha}, \eqref{eq:defv}, the L'H\"ospital rule and the elementary equality $2\sqrt{ab}=a+b-(\sqrt{a}-\sqrt{b})^2$ for $a,b\geq 0$, we have		
		\begin{equation}\label{eq:deltainf}
			\begin{split}
				\delta_{\infty}(\gamma)&=\frac{-n+2\lim\limits_{r\to +\infty}{\alpha(r)}+2k+2}{k+1}\\
				&=\frac{-n+2\left(v_{\infty}+\sqrt{\left(\frac{v_{\infty}+1}{k}\right)\left(k(v_{\infty}-1)+n\right)}\right)+2(k+1)}{k+1}\\
				&=\frac{(k+1)^2\left(\frac{v_{\infty}+1}{k}\right)-\left(\sqrt{k(v_{\infty}-1)+n}-\sqrt{\frac{v_{\infty}+1}{k}}\right)^2}{k+1}\\
				&=\frac{\left(2\left(k+\frac{2(k+1)(v_{\infty}+1)}{k}\right)-n\right)\left(k\sqrt{\frac{v_{\infty}+1}{k}}+\sqrt{k(v_{\infty}-1)+n}\right)}{(k+1)\left((k+2)\sqrt{\frac{v_{\infty}+1}{k}}+\sqrt{k(v_{\infty}-1)+n}\right)}\\
				&=C\left(2\left(k+\frac{2\gamma}{k}+4\right)-n\right),
			\end{split}
		\end{equation}
		where $C$ is a positive constant depending on $n,\,k$ and $\gamma$.
		
		Additionally, from \eqref{ineq:Malpm} and \eqref{eq:deltank}, we obtain			
		\begin{equation}\label{ineq:alpMdinf}
			\delta(r)\geq \delta_{\infty}(\gamma),\,\forall r\geq1.
		\end{equation}
				
		Now, according to the dimension $n$, we consider three cases:
		\begin{itemize}
			\item Case $n>2\left(k+\frac{2\gamma}{k}+4\right)$. By \eqref{eq:deltainf}, $\delta_{\infty}(\gamma)<0$, and by continuity there exists $r_0>1$ such that $\delta(r)<0$ for any $r\geq r_0$. At this point, we have two subcases:
			
			\begin{itemize}
				\item $\lim\limits_{r\to +\infty}{\abs{u(r)}}\in (0,\infty]$. From \eqref{ineq:alpMdinf}, we have $0>\delta(r)\geq \delta_{\infty}(\gamma)$ for any $r\geq r_0$. Then, $1>r^{\delta(r)}\geq r^{\delta_{\infty}(\gamma)}$, $\forall r\geq r_0$ and it follows that 
				\begin{equation*}
					\lim\limits_{r\to
						+\infty}r^{\delta(r)}\in[0,1].
				\end{equation*}
				
				If $\lim\limits_{r\to
					+\infty}r^{\delta(r)}\in(0,1]$, we have a contradiction. Hence, $\lim\limits_{r\to
					+\infty}r^{\delta(r)}=0$ and \eqref{ingrowth} follows
				immediately.
				\item $\lim\limits_{r\to +\infty}{\abs{u(r)}}=0$. Let $R\geq 2r$ and $r\geq 1$. Thus, by the monotony of $u$ and Proposition \ref{prorand2r} there exists $K'>0$ such that:
				\begin{equation}\label{ineq:uRrthm12}
					\begin{split}
						\abs{u(R)-u(r)}&=\abs{u(R)-u(2r)}+\abs{u(2r)-u(r)}\\
						&\geq\abs{u(2r)-u(r)}\geq K'r^{\delta(r)}.
					\end{split}
				\end{equation}
				
				Letting $R\to+\infty$, \eqref{ingrowth} is proved for $r_0=1$.
				
			\end{itemize}
			\item Case $n<2\left(k+\frac{2\gamma}{k}+4\right)$. We have three subcases:
			\begin{itemize}
				\item $\lim\limits_{r\to +\infty}{\abs{u(r)}r^{-\delta(r)}}=0$. From \eqref{rand2r}, we have that $K'=0$, a contradiction.
				\item $\lim\limits_{r\to +\infty}{\abs{u(r)}r^{-\delta(r)}}=L\in(0,+\infty)$. Then, $\forall\epsilon>0$ there exists $r_0\geq 1$ such that $\abs{\abs{u(r)}r^{-\delta(r)}-L}<\epsilon$ for any $r\geq r_0$ and \eqref{ingrowth} is proved for $r_0$.
				\item $\lim\limits_{r\to +\infty}{\abs{u(r)}r^{-\delta(r)}}=+\infty$. Then there exists $r_0\geq 1$ and $K>0$ such that $\abs{u(r)}r^{-\delta(r)}\geq K$ for any $r\geq r_0$  and \eqref{ingrowth} is proved for $r_0$.
			\end{itemize} 	
			\item Case $n=2\left(k+\frac{2\gamma}{k}+4\right)$. Let $r\geq 1$. Then there exists $m\in\natu$ and $1\leq r_1<2$
			such that $r=2^{m-1}r_1$. Thus, by the monotony of $u$ and Proposition \ref{prorand2r}, it follows that
			\begin{equation}\label{ineq:u}
				\begin{split}
					\abs{u(r)}&\geq\abs{u(r)-u(r_1)}-\abs{u(r_1)}\\
					&=\sum_{j=1}^{m-1}\abs{u(2^j r_1)-u(2^{j-1} r_1)}-\abs{u(r_1)}\\
					&\geq \sum_{j=1}^{m-1}{K'(2^{j-1}r_1)^{\delta(2^{j-1}r_1)}}-\abs{u(r_1)}.
				\end{split}
			\end{equation}
			
			By \eqref{ineq:alpMdinf}, we have  
			\begin{equation*}
				\left(2^{j-1}r_1\right)^{\delta(2^{j-1}r_1)}\geq\left(2^{j-1}r_1\right)^{\delta_{\infty}(\gamma)}=1,
			\end{equation*}
			for any $j\in\{1,\ldots,m-1\}$ and together with \eqref{ineq:u}, we get
			\begin{equation}\label{ineq:neq}
					\abs{u(r)}\geq K'(m-1)-\abs{u(r_1)}=\left(\frac{K'}{\log 2}\right)\left(\log r-\log r_1\right) -\abs{u(r_1)},
			\end{equation}
			and \eqref{ingrowth} follows easily.\qedhere
				
		\end{itemize}
	\end{proof}

	\begin{proof}[Proof of Theorem \ref{bogrowth}]
		
		From \eqref{ingrowth} of Theorem \ref{nbgrowth}, it follows that $n\neq2\left(k+\frac{2\gamma}{k}+4\right)$.
		
		Recall that 
		\begin{equation*}
			\delta(r)=\frac{-n+2\alpha(r)+2k+2}{k+1}.
		\end{equation*}

		Let $R\geq 2r$ and $r\geq 1$. From \eqref{ineq:alpMdinf} and \eqref{ineq:uRrthm12}, we have		
		\begin{equation*}
			\abs{u(R)-u(r)}\geq K'r^{\delta(r)}\geq K'r^{\delta_{\infty}(\gamma)}.
		\end{equation*}		
		
		Thus, letting $r\to+\infty$, we conclude that $\delta_{\infty}(\gamma)$ must be negative, which is equivalent to $n>2\left(k+\frac{2\gamma}{k}+4\right)$. This prove $i)$. 
		
		Finally, letting $R\to+\infty$, we have		
		\begin{equation*}
			\abs{u(r)-u_{\infty}}\geq K'r^{\delta(r)},\,\forall r\geq 1,
		\end{equation*} 
		which is $ii)$ with $M=K'$. The proof is complete.
		
	\end{proof}
	
\section*{Appendix}
	We will see that the results obtained in the previous section are optimal.
	
	\begin{example}\label{exam:family}
		Let $\sigma_1,\sigma_2,\tau,\mu\in\real$, $(n\geq 2)\in\natu$ and $k\in\{1,2,\ldots,n\}$. Define the function 
		\begin{equation}\label{eq:wheight}
			w(r)=r^{\sigma_1}\left(1+r^{\sigma_2}\right)^{-\frac{\tau}{\sigma_2}},
		\end{equation}
		and
		\begin{equation}\label{eq:mu}
			\mu=\frac{\sigma_1}{k}+2.
		\end{equation}
		
		We assume that $\sigma_1,\sigma_2,\tau$ and $n$ satisfy the following conditions:
		\begin{itemize}
			\item If $\tau=0$, then
			\begin{equation}\label{ineq:ncond}
				\sigma_1>-k\mbox{ and }\,n>\frac{k(k+1)}{2k+\sigma_1}.
			\end{equation}
			\item If $\tau>0$, then 
			\begin{equation}\label{ineq:ncond1}
				\begin{gathered}
					k(\sigma_2-2)\geq \sigma_1>-k,\,2k+\sigma_1>\tau\mbox{ and }\\
					n>\max\left\lbrace \sqrt{\frac{\tau\sigma_2}{k+1}}+\frac{2k-(k-1)\sigma_1}{k},\,\frac{k(k+1)}{2k+\sigma_1-\tau}\right\rbrace.
				\end{gathered}
			\end{equation}	
		\end{itemize}
		
		Now, let $u_{\beta}$ be a radial function defined by
		\begin{equation*}
			u_{\beta} (r):=\begin{cases}
				\sgn{\beta}\left(1+r^{\mu}\right )^{\frac{\beta}{\mu}} &\mbox{if } \beta
				\neq 0,\forall r\geq 0,\\\\
				\left(\frac{1}{\mu}\right)\log{\left (1+r^{\mu}\right )}&\mbox{if } \beta=0,\forall r\geq 0.
			\end{cases}
		\end{equation*}
		
		On the other hand, let $\beta\in\real$ such that
		\begin{equation}\label{ineq:betaexam}
			\begin{split}
				\beta\geq Q(r)&:=\frac{-n+2\nu+2k+2}{k+1}+\\
				&+\frac{2}{k+1}\begin{cases}
					\frac{\int_1^r{\frac{v(s)}{s}}\,ds}{\log r}&\mbox{if }r\geq 1,\\
					v(r)&\mbox{if }r\in[0,1),
				\end{cases}
			\end{split}
		\end{equation}
		where $\nu:=\sqrt{v(0)^2+\frac{n}{k}v(0)+\frac{n}{k}-1 -v'(1)}$ and $v$ is defined by \eqref{eq:defv}.
		
		Then $u_{\beta}$ is a semi-stable nonconstant radial solution of \eqref{mainequation}
		with $g=g_{\beta}$ defined by
		\begin{itemize}
			\item If $\beta\neq0$,
			\begin{equation*}
				g_{\beta}(s):=c_{n,k}
				\begin{cases}
					\abs{\beta}^k\left(\left(\abs{s}^{\frac{\mu}{\beta}}-1\right)^{\frac{\sigma_2}{\mu}}+1\right)^{\frac{\tau}{\sigma_2}}\times&\mbox{if }s\in I_{\beta\neq0},\\
					\times\left(n+k(\beta-2)+k\left(\mu-\beta\right)\abs{s}^{-\frac{\mu}{\beta}}\right)\abs{s}^{k-\frac{k\mu}{\beta}}\\
					C^1-\mbox{extension}&\mbox{if }s\not\in I_{\beta\neq0},
				\end{cases}
			\end{equation*}
			\item If $\beta=0$,
			\begin{equation*}
				g_{0}(s):=c_{n,k}
				\begin{cases}
					\left(\left(e^{\mu s}-1\right)^{\frac{\sigma_2}{\mu}}+1\right)^{\frac{\tau}{\sigma_2}}\left(n+k(\beta-2)+k\mu e^{-\mu s}\right)e^{-k\mu s}&\mbox{if }s\in I_{0},\\
					C^1-\mbox{extension}&\mbox{if }s\not\in I_{0},
				\end{cases}
			\end{equation*}
		\end{itemize}
		where
		\begin{equation}\label{eq:defIb}
			I_{\beta}:=
			\begin{cases}
				[1,+\infty)&\mbox{if }\beta>0,\\
				[0,+\infty)&\mbox{if }\beta=0,\\
				[-1,0)&\mbox{if }\beta<0.
			\end{cases}
		\end{equation}
	\end{example}
	
	To establish the above result we need the following auxiliary Lemmata.
	\begin{lemma}\label{lemm:theorem_hardygen0}
		Let $\theta\in \real$ and $\rho\in C(\RR^n),\,0\leq V\in C^1\left(\real^n\setminus\{0\}\right)$ be radial functions such that
		\begin{equation}\label{ineq:condV0}
			\theta\left(rV'+(n-2\rho-2)V-\theta V\right)\geq 0,\,\forall r> 0,
		\end{equation}
		and
		\begin{equation}\label{eq:condV1}
			\lim_{r\to 0}{r^{n-2}V}=0.
		\end{equation}
		
		Then
		\begin{equation}\label{ineq:hardygen0}
			\int_0^{\infty}{r^{n-3}V\left(\left(r\eta'+\rho\eta\right)^2-\frac{\theta^2\eta^2}{4}\right) \,dr}\geq 0,
		\end{equation}
		for every radially symmetric function $\eta\in C^1_c(\real^n)$.
	\end{lemma}
	
	\begin{proof}
		Let $\eta\in C_c^{1}(\real^n)$ be a radial function, then
		\begin{equation*}
			\begin{gathered}
				\int_0^{\infty}{r^{n-3}V\left (\theta\eta-t\left(r\eta'+\rho\eta\right)\right )^2\,dr}\geq 0,
			\end{gathered}
		\end{equation*}
		for all $t\in\real$. Extending the above expression, we get the following quadratic inequality for $t$:		
		\begin{equation*}
			\begin{gathered}
				\theta^2\int_0^{\infty}{r^{n-3}V\eta^2\,dr}-t\theta\int_0^{\infty}{r^{n-3}V\left(r(\eta^2)'+2\rho\eta^2\right)\,dr}+\\
				+t^2\int_0^{\infty}{r^{n-3}V\left(r\eta'+\rho\eta\right)^2\,dr}\geq 0.
			\end{gathered}
		\end{equation*}
		
		Integrating by parts and using \eqref{eq:condV1}, we obtain
		\begin{equation*}
			\begin{gathered}
				\theta^2\int_0^{\infty}{r^{n-3}V\eta^2\,dr}+t\theta\int_0^{\infty}{r^{n-3}\eta^2\left(rV'+(n-2\rho-2)V\right)\,dr}+\\
				+t^2\int_0^{\infty}{r^{n-3}V\left(r\eta'+\rho\eta\right)^2\,dr}\geq 0.
			\end{gathered}
		\end{equation*}
		
		Therefore, the above quadratic inequality is equivalent to
		\begin{equation*}
			\begin{split}
				&4\left(\theta^2\int_0^{\infty}{r^{n-3}V\eta^2\,dr}\right)\left(\int_0^{\infty}{r^{n-3}V\left(r\eta'+\rho\eta\right)^2\,dr}\right) \\ &\geq\left(\theta\int_0^{\infty}{r^{n-3}\eta^2\left(rV'+(n-2\rho-2)V\right)\,dr}\right)^2,
			\end{split}
		\end{equation*}
		from \eqref{ineq:condV0}, it follows \eqref{ineq:hardygen0}.
	\end{proof}
	
	We are now ready to establish Example \ref{exam:family}.
	\begin{proof}[Proof of Example \ref{exam:family}]
		
		We claim that $2k+\gamma>0$, where $\gamma$ is given in \eqref{ineq:H1rwrw}. To this end, let $w(r)$ as in \eqref{eq:wheight}. Then, differentiating $\log w(r)=\sigma_1\log r-\frac{\tau}{\sigma_2}\log\left(1+r^{\sigma_2}\right)$ with respect to $r$, we obtain
		\begin{equation*}
			\frac{rw'(r)}{w(r)}=\sigma_1-\tau h_{\sigma_2}(r),
		\end{equation*}
		where
		\begin{equation}\label{eq:hexam}
			h_{\lambda}(r)=\frac{r^{\lambda}}{1+r^{\lambda}}\in[0,1),\,\forall\lambda>0\mbox{ and }r\geq0.
		\end{equation}
		
		It follows that
		\begin{equation}\label{eq:dlhexam}
			\frac{rh'_{\lambda}(r)}{\lambda}=\left(1-h_{\lambda}(r)\right)h_{\lambda}(r)\in\left[0,\frac14\right] ,\,\forall\lambda>0\mbox{ and }r\geq0.
		\end{equation}
		
		From \eqref{eq:defv}, we have
		\begin{equation}\label{eq:vexam}
			v(r)=\frac{k-1+\sigma_1-\tau h_{\sigma_2}(r)}{k+1}=v(0)-\frac{\tau}{k+1}h_{\sigma_2}(r),
		\end{equation}
		with
		\begin{equation*}
			v(0)=\frac{k-1+\sigma_1}{k+1}.
		\end{equation*}
		
		Moreover, we also obtain that
		\begin{equation}\label{eq:dvexam}
			rv'(r)=-\frac{\tau\sigma_2}{k+1}\left(1-h_{\sigma_2}(r)\right)h_{\sigma_2}(r).
		\end{equation}
		
		On the other hand, by \eqref{eq:vexam} we get
		\begin{equation}\label{eq:gamexam}
			\gamma=\lim_{r\to +\infty}{\frac{rw'(r)}{w(r)}}=\lim_{r\to +\infty}{\left\lbrace \sigma_1-\tau h_{\sigma_2}(r)\right\rbrace }=\sigma_1-\tau,
		\end{equation}
		which, by \eqref{ineq:ncond} and \eqref{ineq:ncond1}, implies that 
		\begin{equation}\label{ineq:gamexam}
			2k+\gamma=2k+\sigma_1-\tau>0,
		\end{equation}
		and the claim follows.
		
		Next, we divide the proof into two steps.
		
		\vspace*{1em}\noindent\textbf{Step 1.} For any $\beta\in\real$, $u_{\beta}$ is a $k$-convex solution of \eqref{mainequation} with $g=g_{\beta}$ and $w$ defined by \eqref{eq:wheight}.
		
		A direct calculation gives that
		\begin{equation*}
			u'_{\beta}=\begin{cases}
				\abs{\beta}r^{\mu-1}\left (1+r^{\mu}\right )^{\frac{\beta}{\mu}-1} &\mbox{if } \beta
				\neq 0,\\
				r^{\mu-1}\left (1+r^{\mu}\right )^{-1}&\mbox{if } \beta=0.
			\end{cases}
		\end{equation*}
		
		By \eqref{eq:mu}-\eqref{ineq:ncond1}, it follows that 
		\begin{equation*}
			\mu-1=\frac{\sigma_1}{k}+1>0\Rightarrow u'_{\beta}(0)=0,
		\end{equation*}
		also
		\begin{equation}\label{eq:L2famstable}
			\lambda_{2,\beta}:=\frac{u'_{\beta}}{r}=\begin{cases}
				\abs{\beta}r^{\mu-2}\left (1+r^{\mu}\right )^{\frac{\beta}{\mu}-1} &\mbox{if } \beta
				\neq 0,\\\\
				r^{\mu-2}\left (1+r^{\mu}\right )^{-1}&\mbox{if } \beta=0.
			\end{cases}
		\end{equation}	
		
		Consequently, differentiating $\log\lambda_{2,\beta}$ with respect to $r$ and using \eqref{eq:hexam}, we obtain		
		\begin{equation}\label{eq:rdlamblambexam}
			\frac{r\lambda'_{2,\beta}}{\lambda_{2,\beta}}=\mu-2+\left(\beta-\mu\right)h_{\mu}(r).
		\end{equation}
		
		Fix any  $j\in\{1,2,\ldots,k\}$. By \eqref{eq:radk-Hess}, we have
		\begin{equation}\label{eq:Sjfamstable}
			\begin{split}
				S_j(D^2u_{\beta})&=c_{n,j}\lambda_{2,\beta}^{j-1}\left(n\lambda_{2,\beta}+jr\lambda'_{2,\beta}\right)\\
				&=c_{n,j}\lambda_{2,\beta}^{j}\left(n+j\left((\mu-2)+\left(\beta-\mu\right)h_{\mu}(r)\right)\right)\\
				&=c_{n,j}\lambda_{2,\beta}^{j}((n+j(\mu-2))(1-h_{\mu}(r))+(n+j(\beta-2))h_{\mu}(r)).
			\end{split}
		\end{equation}
		
		Combining \eqref{ineq:vn2k}, \eqref{ineq:ncond}, \eqref{ineq:ncond1} and \eqref{eq:vexam}, we have
		\begin{equation*}
			\begin{split}
				\left(v(0)+\frac{n}{2k}\right)^2-\left(v(r)+\frac{n}{2k}\right)^2&=\left(v(0)-v(r)\right)\left(v(r)+v(0)+\frac{n}{k}\right)\\
				&\geq\left(\frac{2\tau}{k+1}\right)h_{\sigma_2}(r)\left(v(r)+\frac{n}{2k}\right)\\
				&\geq 0,\,\forall r\geq 0.
			\end{split}
		\end{equation*}
		
		By \eqref{eq:hexam} and setting $\lambda=\sigma_2$ in \eqref{eq:dlhexam} together with \eqref{eq:dvexam}, we get
		\begin{equation*}
			\begin{split}
				rv'(r)-v'(1)&=\frac{\tau\sigma_2}{k+1}\left(\frac14-\left(1-h_{\sigma_2}(r)\right)h_{\sigma_2}(r)\right)\geq 0,\,\forall r\geq 0.
			\end{split}
		\end{equation*}
		
		From the last two inequalities and \eqref{ineq:nkvtvtc1}, \eqref{ineq:ncond}, \eqref{ineq:ncond1}, and the fact that $-rv'(r)\geq 0$ for any $r\geq 0$, a straightforward calculation gives
		\begin{equation}\label{ineq:sqrtvexam}
			\begin{split}
				\nu&=\sqrt{v(0)^2+\frac{n}{k}v(0)+\frac{n}{k}-1 -v'(1)}\\
				&=\sqrt{\left(v(0)+\frac{n}{2k}\right)^2 -v'(1)-\left(\frac{n}{2k}-1\right)^2}\\
				&\geq \sqrt{v(r)^2+\frac{n}{k}v(r)+\frac{n}{k}-1-rv'(r)}\geq\abs{v(r)},\,\forall r\geq 0. 
			\end{split}
		\end{equation}
		
		From the previous inequality, \eqref{eq:mu} and \eqref{ineq:betaexam}, we deduce that
		\begin{equation*}
			\begin{split}
				n+j(\mu-2)&=n+\frac{j\sigma_1}{k}>n-j\geq 0,\\
				n+j(\beta-2)&\geq\left(1-\frac{j}{k+1}\right)n+\frac{2\left(v(r)+\abs{v(r)}\right)}{k+1}>0.
			\end{split}
		\end{equation*}
		
		Therefore, from \eqref{eq:Sjfamstable}, $S_j(D^2 u_{\beta})\geq 0$ for any $j\in\{1,2,\ldots,k\}$.	 This show that the functions $u_\beta$ are $k$-convex.
		
		From \eqref{eq:L2famstable} and \eqref{eq:Sjfamstable}, we have
		\begin{equation*}
			\begin{split}
				S_k(D^2u_{\beta})&=c_{n,k}\lambda_{2,\beta}^{k}\left(n+k(\mu-2)+k\left(\beta-\mu\right)\left(\frac{r^{\mu}}{1+r^{\mu}}\right)\right)\\
				&=c_{n,k}\left(\begin{cases}
					\abs{\beta}^k&\mbox{if }\beta
					\neq 0,\\
					1&\mbox{if } \beta=0.
				\end{cases}\right)r^{k(\mu-2)}\left (1+r^{\mu}\right )^{\frac{k\beta}{\mu}-k}\times\\
				&\times\left(n+k(\mu-2)+k\left(\beta-\mu\right)\left(1-\left(1+r^{\mu}\right)^{-1}\right)\right)\\
				&=c_{n,k}\left(\begin{cases}
					\abs{\beta}^k&\mbox{if }\beta
					\neq 0,\\
					1&\mbox{if } \beta=0.
				\end{cases}\right)r^{k(\mu-2)}\left(1+r^{\sigma_2}\right)^{\frac{-\tau}{\sigma_2}}\times\\
				&\times\left(1+r^{\sigma_2}\right)^{\frac{\tau}{\sigma_2}} \left(\left(n+k(\beta-2)\right)\left (1+r^{\mu}\right )^{\frac{k\beta}{\mu}-k}\right.+\\
				&-\left.k\left(\beta-\mu\right)\left(1+r^{\mu}\right)^{\frac{k\beta}{\mu}-(k+1)}\right).
			\end{split}
		\end{equation*}
		
		On the other hand, it is easy to see that
		\begin{equation*}
			1+r^{\mu}=\begin{cases}
				\abs{u_{\beta}(r)}^{\frac{\mu}{\beta}} &\mbox{if } \beta
				\neq 0,\forall r\geq 0,\vspace*{5pt}\\
				e^{\mu u_{0}(r)}&\mbox{if } \beta=0,\forall r\geq 0,
			\end{cases}\mbox{ and }
			r^{\sigma_2}=\begin{cases}
				\left(\abs{u_{\beta}(r)}^{\frac{\mu}{\beta}}-1\right)^{\frac{\sigma_2}{\mu}} &\mbox{if } \beta
				\neq 0,\forall r\geq 0,\vspace*{5pt}\\
				\left(e^{\mu u_{0}(r)}-1\right)^{\frac{\sigma_2}{\mu}}&\mbox{if } \beta=0,\forall r\geq 0.
			\end{cases}
		\end{equation*}
		
		From this, \eqref{eq:wheight} and \eqref{eq:mu}, we have 
		\begin{equation*}
			\begin{split}
				S_k(D^2u_{\beta})&=w(r)c_{n,k}\begin{cases}
					\abs{\beta}^{k}\left(\left(\abs{u_{\beta}}^{\frac{\mu}{\beta}}-1\right)^{\frac{\sigma_2}{\mu}}+1\right)^{\frac{\tau}{\sigma_2}}\times&\mbox{if }\beta\neq 0,\\
					\times\left(n+k(\beta-2)+k\left(\mu-\beta\right)\abs{u_{\beta}}^{-\frac{\mu}{\beta}}\right)\abs{u_{\beta}}^{k\left(\frac{\beta-\mu}{\beta}\right)}\\
					\left(\left(e^{\mu u_{0}}-1\right)^{\frac{\sigma_2}{\mu}}+1\right)^{\frac{\tau}{\sigma_2}}\left(n-2k+k\mu e^{-\mu u_{0}}\right)e^{-k\mu u_{0}}&\mbox{if }\beta=0.
				\end{cases}
			\end{split}
		\end{equation*}
		
		Now for every $s\geq 0$, let us define the function
		\begin{equation}\label{eq:hbetexam}
			f_{\beta}(s):=\left(q_{\beta}(s)\right)^{\frac{\tau}{\sigma_2}}\begin{cases}
				\abs{\beta}^{k}
				\left(n+k(\beta-2)+k\left(\mu-\beta\right)s^{-\frac{\mu}{\beta}}\right )s^{k\left(\frac{\beta-\mu}{\beta}\right)}&\mbox{if }\beta\neq 0,\\
				\left(n-2k+k\mu e^{-\mu s}\right)e^{-k\mu s}&\mbox{if }\beta=0,
			\end{cases}
		\end{equation}
		with
		\begin{equation}\label{eq:qbetexam}
			q_{\beta}(s):=\begin{cases}
				\left(s^{\frac{\mu}{\beta}}-1\right)^{\frac{\sigma_2}{\mu}}+1&\mbox{if }\beta\neq 0,\\
				\left(e^{\mu s}-1\right)^{\frac{\sigma_2}{\mu}}+1&\mbox{if }\beta=0,
			\end{cases}
		\end{equation}
		which lead us to
		\begin{equation*}
			S_k(D^2u_{\beta})=w(r)c_{n,k}\begin{cases}
				f_{\beta}(\abs{u_{\beta}})&\mbox{if }\beta\neq0,\\
				f_{0}(u_{0})&\mbox{if }\beta=0.
			\end{cases}
		\end{equation*}
		
		Let $y\in I_{\beta}$ and $s=\abs{y}$. Using \eqref{eq:defIb}, we see that:
		\begin{itemize}
			\item If $\beta>0$, then $s\geq 1$.
			\item If $\beta=0$, then $s\geq 0$.
			\item If $\beta<0$, then $s\in(0,1]$.
		\end{itemize}
		
		Therefore, to study the differentiability of $f_{\beta}$, we must to consider the points $s=0$ for $\beta\leq 0$ and $s=1$ for $\beta\neq0$.
		
		Since $q_{\beta\neq0}(1)=q_0(0)=1$, from \eqref{eq:hbetexam} and \eqref{eq:qbetexam}, we have 
		\begin{equation}\label{eq:cont0}
			\lim_{s\to 1}{\left. \frac{f_{\beta}(s)}{\abs{\beta}^k}\right|_{\beta\neq 0} }=\lim_{s\to 0}{f_{0}(s)}=n+k(\mu-2),\,\forall\tau\geq0.
		\end{equation}		
		
		Now, differentiating $\log f_{\beta}(s)$ with respect to $s$ and using \eqref{eq:hbetexam}, we obtain 
		\begin{equation}\label{eq:dhbetexam}
			\frac{f'_{\beta}(s)}{f_{\beta}(s)}=\begin{cases}
				0&\mbox{if }\tau=0,\\
				\frac{\tau}{\sigma_2}\frac{q'_{\beta}(s)}{q_{\beta}(s)}&\mbox{if }\tau>0.
			\end{cases}+\begin{cases}
				\frac{k(\beta-\mu)s^{-1}}{\beta}\left(1+\frac{\mu s^{-\frac{\mu}{\beta}}}{n+k(\beta-2)+k\left(\mu-\beta\right)s^{-\frac{\mu}{\beta}}}\right)&\mbox{if }\beta\neq 0,\vspace*{5pt}\\
				-k\mu\left(1+\frac{\mu e^{-\mu s}}{n-2k+k\mu e^{-\mu s}}\right) &\mbox{if }\beta=0.
			\end{cases}
		\end{equation}
		
		If $\tau=0$, then from \eqref{eq:cont0} and \eqref{eq:dhbetexam}, we get		
		\begin{equation}\label{eq:cont02}
			\lim_{s\to 1}{\left.\frac{f'_{\beta}(s)}{\abs{\beta}^k}\right|_{\beta\neq 0}}=\frac{k(\beta-\mu)((k+1)\mu+n-2k)}{\beta},\,\lim_{s\to 0}{f'_{0}(s)}=-k\mu((k+1)\mu+n-2k),
		\end{equation}		
		and if $\tau>0$, from \eqref{eq:qbetexam}, we have
		\begin{equation}\label{eq:dqbetexam}
			q'_{\beta}(s)=\sigma_2\begin{cases}
				\beta^{-1}\left(s^{\frac{\mu}{\beta}}-1\right)^{\frac{\sigma_2}{\mu}-1}s^{\frac{\mu}{\beta}-1}&\mbox{if }\beta\neq0,\\
				\left(e^{\mu s}-1\right)^{\frac{\sigma_2}{\mu}}e^{\mu s}&\mbox{if }\beta=0.
			\end{cases}
		\end{equation}
		
		From \eqref{eq:mu} and \eqref{ineq:ncond1}, it follows that $\mu\leq\sigma_2$ and then 		
		\begin{equation}\label{eq:cont1}
			\begin{split}
				\lim_{s\to 1}{q'_{\beta\neq0}(s)}=\begin{cases}
					\frac{\sigma_2}{\beta}&\mbox{if }\mu=\sigma_2,\\
					0&\mbox{if }\mu<\sigma_2.
				\end{cases},\,\lim_{s\to 0}{q'_{0}(s)}=\begin{cases}
					\sigma_2&\mbox{if }\mu=\sigma_2,\\
					0&\mbox{if }\mu<\sigma_2.
				\end{cases}
			\end{split}
		\end{equation}
		
		Therefore, concerning the cases $s=0$ for $\beta= 0$ and $s=1$ for $\beta\neq0$, we have  		
		\begin{equation}\label{eq:cont2}
			\begin{split}
				\lim_{s\to 1}{\left.\frac{f'_{\beta}(s)}{\abs{\beta}^k}\right|_{\beta\neq 0}}&=\left(\frac{\tau}{\sigma_2}\lim_{s\to 1}{q'_{\beta\neq0}(s)}\right)(n+k(\mu-2))+\frac{k(\beta-\mu)((k+1)\mu+n-2k)}{\beta},\\
				\lim_{s\to 0}{f'_{0}(s)}&=\left(\frac{\tau}{\sigma_2}\lim_{s\to 0}{q'_{0}(s)}\right)(n+k(\mu-2))-k\mu((k+1)\mu+n-2k).
			\end{split}
		\end{equation}
		
		For the case when
		$\beta<0$ and $\tau\geq0$, we rewrite $q_{\beta}(s)$ to get
		\begin{equation}\label{eq:cont3}
			q_{\beta}(s)=\left(\left(1-s^{\frac{-\mu}{\beta}}\right)^{\frac{\sigma_2}{\mu}}+s^{\frac{-\sigma_2}{\beta}}\right)s^{\frac{\sigma_2}{\beta}} \mbox{ and }
			\frac{q'_{\beta}(s)}{q_{\beta}(s)}=\frac{\sigma_2}{\beta}\left(\frac{\left(1-s^{\frac{-\mu}{\beta}}\right)^{\frac{\sigma_2}{\mu}-1}}{\left(1-s^{\frac{-\mu}{\beta}}\right)^{\frac{\sigma_2}{\mu}}+s^{\frac{-\sigma_2}{\beta}}}\right)s^{-1}.
		\end{equation}
		
		From \eqref{eq:mu} and \eqref{ineq:gamexam}, it follows that $k\mu-\tau=2k+\sigma_1-\tau>0$. Then
		\begin{equation}\label{ineq:mtb}
			\frac{\tau-k\mu}{\beta}=-\frac{k\mu-\tau}{\beta}>0.
		\end{equation}
		
		Thus, from \eqref{eq:hbetexam}, \eqref{eq:dhbetexam} and \eqref{eq:cont3}, for $\beta\neq0$, we get
		\begin{equation*}
			\begin{split}
				\frac{f_{\beta}(s)}{\abs{\beta}^{k}}=\left(n+k(\beta-2)+o(1)\right)s^{k+\frac{\tau-k\mu}{\beta}},\,\frac{f'_{\beta}(s)}{f_{\beta}(s)}=\left(\frac{\tau+k(\beta-\mu)}{\beta}+o(1)\right)s^{-1},
			\end{split}
		\end{equation*}
		as $s\to0$.
		
		Therefore, from \eqref{ineq:mtb}, we obtain
		\begin{equation}\label{eq:cont4}
			\lim_{s\to0}{\left.\frac{ f_{\beta}(s)}{\abs{\beta}^k}\right|_{\beta<0}}=\left(n+k(\beta-2)\right)\lim_{s\to0}{s^{k+\frac{\tau-k\mu}{\beta}}}=0,
		\end{equation}
		and
		\begin{equation}\label{eq:cont5}
			\lim_{s\to0}{\left.\frac{ f'_{\beta}(s)}{\abs{\beta}^k}\right|_{\beta<0}}=\frac{\left(n+k(\beta-2)\right)\left(\tau+k(\beta-\mu)\right)}{\beta}\lim_{s\to0}{s^{k-1+\frac{\tau-k\mu}{\beta}}}=0.
		\end{equation}
		
		Collecting \eqref{eq:cont0}, \eqref{eq:cont02}, \eqref{eq:cont2}, \eqref{eq:cont4} and \eqref{eq:cont5}, there exists a $C^1$-extension for $f_{\beta}(s)$ when $s\not\in I_{\beta}$, and we have that $u_{\beta}$ is a radial solution of \eqref{mainequation} with $g=g_{\beta}$, where
		\begin{equation}
			g_{\beta}(s)=c_{n,k}\begin{cases}
				\begin{cases}
					f_{\beta}(\abs{s})&\mbox{if }s\in I_{\beta\neq 0},\\
					C^1\mbox{-extension}&\mbox{if }s\not\in I_{\beta\neq 0},
				\end{cases}\\
				\begin{cases}
					f_0(s)&\mbox{if }s\in I_0,\\
					C^1\mbox{-extension}&\mbox{if }s\not\in I_{0}.
				\end{cases}
			\end{cases}
		\end{equation}
		
		\vspace*{1em}\noindent\textbf{Step 2.} For $\beta$ satisfying the inequality \eqref{ineq:betaexam}, $u_{\beta}$ is a semistable solution of \eqref{mainequation}.
		
		Next, for suitable $\rho,\theta$ and $V$ we verify the hypotheses of Lemma \ref{lemm:theorem_hardygen0}. To this end, let $V=r^2\lambda_{2,\beta}^{k+1}$. Differentiating $\log V$ with respect to $r$, from \eqref{eq:hexam}, \eqref{eq:L2famstable} and \eqref{eq:rdlamblambexam}, we have
		\begin{equation}\label{eq:rVVexam}
			\frac{rV'}{V}=2+(k+1)\left(\mu-2+\left(\beta-\mu\right)h_{\mu}(r)\right).
		\end{equation} 
		
		Note that in both cases $\tau=0$ or $\tau>0$, we have that $n>k(k+1)/(2k+\sigma_1)$ and $k+\sigma_1>0$ by \eqref{ineq:ncond} and \eqref{ineq:ncond1}. Now, since
		\begin{equation*}
			\frac{k(k+1)}{2k+\sigma_1}+\frac{(k+1)\sigma_1}{k}=\frac{(k+1)\left(k+\sigma_1\right)^2}{k(2k+\sigma_1)}>0,
		\end{equation*}
		it follows that $n>-(k+1)(\sigma_1/k)$. Combining this with \eqref{eq:mu} and \eqref{eq:L2famstable}, we have		
		\begin{equation}\label{eq:cond32}
			\lim_{r\to 0}{r^{n-2}V}=C_{\beta}\lim_{r\to 0}{r^{n+\frac{(k+1)\sigma_1}{k}}\left (1+r^{\mu}\right )^{\frac{(k+1)(\beta-\mu)}{\mu}}}=0,
		\end{equation}
		with $C_{\beta}=\abs{\beta}^{k+1}$ for $\beta\neq 0$ and $C_{\beta}=1$ for $\beta=0$.
		
		Now, consider $\rho=v$ and $\theta=2\nu$. From \eqref{eq:rVVexam}, we have
		\begin{equation*}
			\begin{split}
				n-2\rho-\theta-2+\frac{rV'}{V}&=n-2\left(v+\nu\right)
				+(k+1)\left(\mu-2+\left(\beta-\mu\right)h_{\mu}(r)\right).
			\end{split}
		\end{equation*}
		
		We have from \eqref{ineq:betaexam} that
		\begin{equation*}
			n=(k+1)(2-Q(r))+2\nu+2\begin{cases}
				\frac{\int_1^r{\frac{v(s)ds}{s}}}{\log r}&\mbox{if }r\geq 1,\\
				v(r)&\mbox{if }r\in[0,1).
			\end{cases}
		\end{equation*}
		
		From this we obtain
		\begin{equation*}
			\begin{split}
				n-2\rho-\theta-2+\frac{rV'}{V}&=(k+1)\left(\mu-Q(r)+\left(\beta-\mu\right)h_{\mu}(r)\right)+\\
				&+2\begin{cases}
					\frac{\int_1^r{\frac{v(s)ds}{s}}}{\log r}-v(r)&\mbox{if }r\geq 1,\\
					0&\mbox{if }r\in[0,1)
				\end{cases}\\
				&=(k+1)\left(\mu-Q(r)\right)\left(1-h_{\mu}(r)\right) +\\
				&+(k+1)\left(\beta-Q(r)\right)h_{\mu}(r)+2\begin{cases}
					\frac{\int_1^r{\frac{v(s)ds}{s}}}{\log r}-v(r)&\mbox{if }r\geq 1,\\
					0&\mbox{if }r\in[0,1).
				\end{cases}
			\end{split}
		\end{equation*}
		
		Consider the function
		\begin{equation*}
		A(r)=2\begin{cases}
			\frac{\int_1^r{\frac{v(s)ds}{s}}}{\log r}-v(r)&\mbox{if }r\geq 1,\\
			0&\mbox{if }r\in[0,1).
		\end{cases}
		\end{equation*}
		
		Note that $\partial_r\left(\int_1^r{\frac{v(s)ds}{s}}-v(r)\log r\right)=-v'(r)\log r\geq 0$ for any $r\geq1$ by \eqref{eq:dvexam}. Then $\int_1^r{\frac{v(s)ds}{s}}-v(r)\log r\geq0$ for every $r\geq1$. Therefore $A(r)\geq 0$ for every $r\geq 0$. Now since, $\beta\geq Q(r)$ and $h_{\mu}(r)\in[0,1)$, it follows that
		\begin{equation}\label{ineq:essexam}
			n-2\rho-\theta-2+\frac{rV'}{V}\geq 0\mbox{ if }(k+1)\left(\mu-Q(r)\right)\left(1-h_{\mu}(r)\right)+A(r)\geq 0.
		\end{equation}
		
		We claim that 
		\begin{equation}\label{eq:nonnega}
			(k+1)\left(\mu-Q(r)\right)\left(1-h_{\mu}(r)\right)+A(r)\geq 0.
		\end{equation}
		
		Using \eqref{eq:mu}, \eqref{eq:vexam} and \eqref{eq:dvexam}, we can easily see that		
		\begin{equation}\label{eq:vexam2}
			\begin{split}
				\mu=\frac{(k+1)(v(0)+1)}{k},\,v(0)-v(r)=\frac{\tau}{k+1}h_{\sigma_2}(r)\mbox{ and }v'(1)=-\frac{\tau\sigma_2}{4(k+1)}.
			\end{split}
		\end{equation}
		
		Let us now consider the functions:
		\begin{equation}\label{eq:Iexam}
			I_{\pm}:=I_{\pm}(r)=\frac{(k^2+1)(v(0)+1)}{k}+n-2k+\frac{2\tau}{k+1}h_{\sigma_2}(r)
			\pm2\nu.
		\end{equation}
		
		Since $n>k/(v(0)+1)$ by \eqref{ineq:nkvtvtc1}, we have		
		\begin{equation}\label{ineq:k21v01n2k}
			\begin{split}
				\frac{(k^2+1)(v(0)+1)}{k}+n-2k&>\left(\frac{k}{v(0)+1}\right)\left(\frac{(k^2+1)(v(0)+1)^2}{k^2}-2v(0)-1\right)\\
				&=\frac{kv(0)^2}{v(0)+1}+\frac{v(0)+1}{k}\\
				&>0.
			\end{split}
		\end{equation}
		
		From this and the fact that $\nu\geq 0$ by \eqref{ineq:sqrtvexam}, it follows that $I_{+}>0$.
		
		On the other hand, from \eqref{ineq:betaexam}, \eqref{eq:vexam2} and \eqref{eq:Iexam}, we have
		\begin{equation}\label{eq:mubetexam}
			\begin{split}
				(k+1)(\mu-Q(r))&=\frac{(k+1)^2(v(0)+1)}{k}+n-2\nu-2k-2+2\begin{cases}
					-\frac{\int_1^r{\frac{v(s)ds}{s}}}{\log r}&\mbox{if }r\geq 1,\\
					-v(r)&\mbox{if }r\in[0,1)
				\end{cases}\\
				&=\frac{(k^2+1)(v(0)+1)}{k}+n-2k+2(v(0)+1)-2(v(r)+1)-2\nu-A(r)\\
				&=I_{-}-A(r).
			\end{split}
		\end{equation}
		
		It follows that
		\begin{equation*}
			\begin{split}
			(k+1)\left(\mu-Q(r)\right)\left(1-h_{\mu}(r)\right)+A(r)&=(I_{-}-A(r))\left(1-h_{\mu}(r)\right)+A(r)\\
			&=\left(1-h_{\mu}(r)\right)I_{-}+A(r)h_{\mu}(r)\\
			&\geq\left(1-h_{\mu}(r)\right)I_{-}.
			\end{split}
		\end{equation*}
		
		In the following, we will prove that $I_{-}I_{+}\geq 0$, which implies that $I_{-}\geq 0$ since $I_{+}>0$. To do that, we begin by noting that
		\begin{equation*}
			\nu^2=v(0)^2+\frac{n}{k}v(0)+\frac{n}{k}-1-v'(1)=\left(v(0)+1\right)\left(v(0)+\frac{n}{k}-1\right)+\frac{\tau\sigma_2}{4(k+1)},
		\end{equation*} 
		and
		\begin{equation}\label{eq:mubet2exam1}
			I_{-}I_{+}=B(r)^2-4\nu^2,
		\end{equation}
		where $B(r):=\frac{(k^2+1)(v(0)+1)}{k}+n-2k+\frac{2\tau}{k+1}h_{\sigma_2}(r)$.
		Using that
		\begin{equation*}
				\left(\frac{(k^2-1)(v(0)+1)}{k}+n-2k\right)^2=\left(B(r)-\frac{2\tau}{k+1}\,h_{\sigma_2}(r)\right)^2-4\left(\nu^2-\frac{\tau\sigma_2}{4(k+1)}\right),
		\end{equation*}
		together with \eqref{ineq:k21v01n2k} and \eqref{eq:mubet2exam1} we then have
		\begin{equation}\label{ineq:ineqIexam}
			\begin{split}
				I_{-}I_{+}-\left(\frac{(k^2-1)(v(0)+1)}{k}+n-2k\right)^2&=\frac{4\tau}{k+1}h_{\sigma_2}(r)\left(\frac{(k^2+1)(v(0)+1)}{k}+n-2k\right)+\\
				&+\frac{4\tau^2}{(k+1)^2}h^2_{\sigma_2}(r)-\frac{\tau\sigma_2}{k+1}\\
				&\geq -\frac{\tau\sigma_2}{k+1}.
			\end{split}
		\end{equation}
		
		If $\tau=0$, we are done. So assume that $\tau>0$, then		
		\begin{equation*}
			\begin{split}
				\left(\frac{(k^2-1)(v(0)+1)}{k}+n-2k\right)^2-\frac{\tau\sigma_2}{k+1}&=\left(\frac{(k^2-1)(2k+\sigma_1)}{k(k+1)}+n-2k-\sqrt{\frac{\tau\sigma_2}{k+1}}\right)\times\\
				&\times\left(\frac{(k^2-1)(2k+\sigma_1)}{k(k+1)}+n-2k+\sqrt{\frac{\tau\sigma_2}{k+1}}\right)\\&\geq 0,
			\end{split}
		\end{equation*}	
		by \eqref{ineq:ncond1}. From this and \eqref{ineq:ineqIexam}, $I_{-}I_{+}\geq 0$, which proves the claim \eqref{eq:nonnega}. 
		
		Thus, according to \eqref{ineq:essexam} and \eqref{eq:nonnega}, we have 		
		\begin{equation*}
			\begin{split}
				\theta\left(rV'+(n-2\rho-2)V-\theta V\right)&=\theta V\left(n-2\rho-\theta-2+\frac{rV'}{V}\right)\geq 0. 
			\end{split}
		\end{equation*}
		
		Hence, from the previous inequality and from \eqref{eq:cond32}, we conclude that the conditions \eqref{ineq:condV0} and \eqref{eq:condV1} of Lemma \ref{lemm:theorem_hardygen0} are satisfied. Recall that $V=r^2\lambda_{2,\beta}^{k+1},\, \rho=v$ and $\theta=2\nu$, then one can see from \eqref{ineq:hardygen0} and \eqref{ineq:sqrtvexam} that		
		\begin{equation*}
			\begin{split}
				\int_0^{\infty}{r^{n-1}\lambda_{2,\beta}^{k+1}\left(r\eta'+v\eta\right)^2}dr&\geq\int_0^{\infty}{r^{n-1}\lambda_{2,\beta}^{k+1}\left(\frac{\theta^2\eta^2}{4}\right) \,dr}\\
				&\geq\int_0^{\infty}{r^{n-1}\lambda_{2,\beta}^{k+1}\left(v^2+\frac{n}{k}v+\frac{n}{k}-1-rv'\right)\eta^2 \,dr},
			\end{split}
		\end{equation*}
		for every radially symmetric function $\eta\in C^1_c(\real^n)$.
		
		Therefore, from Lemma \ref{essential0} and Corollary 1.8 of \cite{NaSa21}, we have that $u_\beta$ is a semi-stable solution of \eqref{mainequation}. Although only the case of a ball appear in Corollary 1.8, one can see that the proof of this lemma can be adapted without difficulties to the $\RR^n$ space. 
	\end{proof}
	
	\begin{proof}[Proof of Remark \ref{rmk1}]
	Let $u_{\beta}$ as in Remark \ref{rmk1}
	a semi-stable solution of \eqref{mainequation}. From Theorem \ref{nbgrowth}, there exist $M>0$ and $r_0\geq 1$, such that for any $r\geq r_0$, we have
		
		\begin{itemize}
			\item If $n\neq 2\left(k+\frac{2\sigma_1}{k}+4\right)$,
				\begin{itemize}
				\item If $\beta\neq 0$, then
				\begin{equation}\label{exam:bn00}
				C_{\beta}r^{\beta-\delta(r)}\geq\frac{\abs{u_{\beta}(r)}}{r^{\delta(r)}}\geq M,
				\end{equation}
				where $C_{\beta}=2^{\beta/\mu}$ if $\beta>0$ and $C_{\beta}=1$ if $\beta<0$, 
				\item If $\beta=0$, then
				\begin{equation}\label{exam:b00}
				\left(\frac{\log 2}{\mu}+\log r\right) r^{-\delta(r)}\geq\frac{\abs{u_{0}(r)}}{r^{\delta(r)}}\geq M,
				\end{equation}
				\end{itemize}			
			\item If $n=2\left(k+\frac{2\sigma_1}{k}+4\right)$ and $\beta\neq 0$, then
				\begin{equation}\label{exam:bn01}
				C_{\beta}\frac{r^{\beta}}{\log r}\geq\frac{\abs{u_{\beta}(r)}}{\log r}\geq M,
				\end{equation}
				where $C_{\beta}$ is as in \eqref{exam:bn00}. 	
		\end{itemize}

Now, let $\tau=0$ in \eqref{eq:wheight}, then $w(r)=r^{\sigma_1}$. Thus from \eqref{eq:deltank} and \eqref{ineq:betaexam} we get
\begin{equation*}
Q(r)=\delta(r)=\delta_{\infty}(\sigma_1),\,\forall r\geq 1.  
\end{equation*}	

So, if $\delta_{\infty}(\sigma_1)>\beta$, we have
\begin{itemize}
	\item If $n\neq2\left(k+\frac{2\sigma_1}{k}+4\right)$, from \eqref{exam:bn00} and \eqref{exam:b00}, letting $r\to+\infty$, we obtain a contradiction.
	\item If $n=2\left(k+\frac{2\sigma_1}{k}+4\right)$, from \eqref{eq:deltainf}, we have that $\delta_{\infty}(\sigma_1)=0$ and from \eqref{exam:bn01}, letting $r\to+\infty$, we obtain a contradiction. 
\end{itemize}		
		
		Hence, $\beta\geq \delta_{\infty}(\sigma_1)=\delta(r)=Q(r),\,\forall r\geq 1$.
		
		In addition, it is easy to check that for any $r\geq 1$, we have that
		
		\begin{itemize}
	\item If $n\neq 2\left(k+\frac{2\sigma_1}{k}+4\right)$,
	\begin{equation*}
	\frac{\abs{u_{\beta}(r)}}{r^{\delta_{\infty}(\sigma_1)}}\geq\begin{cases}
	2^{\frac{\beta}{\mu}}&\mbox{if }\beta<0,\\
	\frac{1}{\mu}\log 2&\mbox{if }\beta=0,\\
	1&\mbox{if }\beta>0.
	\end{cases}
	\end{equation*}
	\item If $n=2\left(k+\frac{2\sigma_1}{k}+4\right)$,
	\begin{equation*}
	\frac{\abs{u_{\beta}(r)}}{\log r}\geq\begin{cases}
	\beta e&\mbox{if }\beta>0,\\
	1&\mbox{if }\beta=0.
	\end{cases}
	\end{equation*}
\end{itemize}		
		
Therefore, we conclude that Theorem \ref{nbgrowth} is sharp for $\gamma=\sigma_1>-k$ and $w(r)=r^{\sigma_1}$ where
\begin{equation*}
M=\begin{cases}
	2^{\frac{\beta}{\mu}}&\mbox{if }\beta<0,\\
	\frac{\log 2}{\mu}&\mbox{if }\beta=0,\\
	\min\{\beta e,1\}&\mbox{if }\beta>0,
	\end{cases}
	\end{equation*}
and $r_0=1$.
\end{proof}
	
	\begin{proof}[Proof of Remark \ref{rmk2}]
	Let $\beta<0$ such that \eqref{ineq:betaexam} holds. From Example \ref{exam:family} $u_{\beta}(r)=-(1+r^\mu)^\frac{\beta}{\mu}$ is a semi-stable bounded solution of \eqref{mainequation}. 
		Now note that from \eqref{ineq:vn2k}, \eqref{eq:vexam}, \eqref{eq:dvexam} and the fact that, for $r\geq1$, $h_{\sigma_2}(r)\in[1/2,1)$, $h'_{\sigma_2}(r)>0$, we have
		\begin{equation*}
			\begin{split}
				\partial_r\left(\sqrt{v(r)^2+\frac{n}{k}v(r)+\frac{n}{k}-1-rv'(r)}\right)&=\left(\frac{-\tau}{k+1}\right)h'_{\sigma_2}(r)\times\\
				&\times\left(\frac{v(r)+\frac{n}{2k}+\sigma_2\left(2h_{\sigma_2}(r)-1\right)}{\sqrt{v(r)^2+\frac{n}{k}v(r)+\frac{n}{k}-1-rv'(r)}}\right)\\
				&<0.
			\end{split}
		\end{equation*}
		
		From this and \eqref{ineq:sqrtvexam}, we obtain
		\begin{equation*}
			\begin{split}
				\nu=\sqrt{v(0)^2+\frac{n}{k}v(0)+\frac{n}{k}-1 -v'(1)}&>\sqrt{v(1)^2+\frac{n}{k}v(1)+\frac{n}{k}-1-v'(1)}\\
				&\geq\frac{\int_1^r{\frac{\sqrt{v(s)^2+\frac{n}{k}v(s)+\frac{n}{k}-1 -sv'(s)}}{s}ds}}{\log r},\,\forall r\geq 1.
			\end{split}
		\end{equation*}
		
		Hence, from \eqref{ineq:alpMdinf}, \eqref{ineq:betaexam} and \eqref{eq:gamexam}, we have
		\begin{equation}\label{ineq:ubrmk}
			\begin{split}
				Q(r)&=\frac{-n+2\nu+2k+2}{k+1}+\left(\frac{2}{k+1}\right)\frac{\int_1^r{\frac{v(s)ds}{s}}}{\log r}\\
				&\geq\frac{-n+2\alpha(r)+2k+2}{k+1}=\delta(r) ,\,\forall r\geq 1.
			\end{split}
		\end{equation}
		
		From this, we have that $0>\beta\geq \delta(r) ,\,\forall r\geq 1$ and item $i)$ follows from \eqref{eq:deltainf}.
		
		Since $\beta-\delta(r)\geq 0$ for any $r\geq 1,\, u_{\beta,\infty}=\lim_{r\to +\infty}{u_{\beta}(r)}=0$ and $\left(1+r^{-\mu}\right)^{\frac{\beta}{\mu}}\geq 2^{\frac{\beta}{\mu}}$, from \eqref{ineq:ubrmk}, we finally obtain
		\begin{equation*}
			\frac{\abs{u_{\beta}(r)-u_{\beta,\infty}}}{r^{\delta(r)}}\geq \frac{\abs{u_{\beta}(r)}-\abs{u_{\beta,\infty}}}{r^{\delta(r)}}=(1+r^{-\mu})^\frac{\beta}{\mu}r^{\beta-\delta(r)}\geq 2^{\frac{\beta}{\mu}},
		\end{equation*}
		and item $ii)$ follows. 
	\end{proof}
\section*{Acknowledgements}
M. Navarro was supported by XUNTA de Galicia under Grant Axudas \'{a} etapa de formaci\'{o}n posdoutoral 2017 and partially supported by AEI of Spain under Grant MTM2016-75140-P and co-financed by European Community fund FEDER and XUNTA de Galicia under grants GRC2015/004, R2016/022 and ED431C2019/02. J. S\'{a}nchez was partially supported by Fondecyt grant 1221928.
	\bibliography{k-Hessianbib}
	\bibliographystyle{plain}
\end{document}